\documentclass[authoryear]{elsarticle}

\journal{arXiv.org}

\usepackage{bbm,bm,latexsym,amsfonts,natbib} 
\newtheorem{Thm}{Theorem}
\newcommand{\startthm} {\begin{Thm}}
\newcommand{\stopthm} {\end{Thm}}

\newtheorem{Cor}[Thm]{Corollary}
\newcommand{\startcor} {\begin{Cor}}
\newcommand{\stopcor} {\end{Cor}}

\newtheorem{Lem}[Thm]{Lemma}
\newcommand{\startlem} {\begin{Lem}}
\newcommand{\stoplem} {\end{Lem}}

\newcommand{\startpf}[1] {\vspace{5mm} \noindent \textbf{Proof #1}  }
\newcommand{\stoppf} {\hfill $\Box$ \\ \vspace{\baselineskip} }

\newcommand{\myitem}[1]{\noindent #1}
\newcommand{\strutfive}{\rule[-2.5mm]{0cm}{5mm}}
\newcommand{\pialphabytwo}{ {{\scriptstyle{\pi \alpha \over 2}}}}
\newcommand{\tanconst} {\tan \pialphabytwo }
\newcommand{\reals} {\mathbb{R}}
\newcommand{\Stable}[3] {\ensuremath{{\mathbf{S}\left(#1, #2 #3\right)}}}
\newcommand{\gt}{\widetilde{g}}

\newcommand{\sign} {{\mathrm{sign} \,}}
\newcommand{\sphere}{\mathbb{S}}
\newcommand{\bolds} {\mathbf{ s}}
\newcommand{\covariation}[3]{{[ #1 , #2 ]_#3}}
\newcommand{\signedpower}[2]{#1^{<#2>}}                   
\newcommand{\sgnpower}[1]{{^{<#1>}}}
\newcommand{\parenspipbytwo}{ {\scriptstyle{ (\frac{\pi p} 2) }}}
\newcommand{\comment}[1]{{\index{zzz: comments}}}
\newcommand{\indicator}[1]{ \mathbbm{1}_#1 }  
\newcommand{\eqd}{{\stackrel{d}{=}}}
\newcommand{\unfinished}{{(unfinished)\index{zzz: unfinished} }}
\newcommand{\citetwo}[3]{{\cite{#1}}}
\newcommand{\mycite}[2]{{\cite{#1}}}

\begin{document}

\begin{frontmatter}

\title{Truncated fractional moments of stable laws}

\author{John P. Nolan\fnref{myfootnote}}
\fntext[myfootnote]{The author was supported by an agreement with
   Cornell University, Operations Research \& Information Engineering under W911NF-12-1-0385
   from the Army Research Development and Engineering Command.}
\address{American University }
\ead{jpnolan@american.edu}

\begin{abstract}
Expressions are given for the truncated fractional moments $E X_+^p$ of a general stable
law. These involve families of special functions that arose out of the study of
multivariate stable densities and probabilities.
As a particular case, an expression is given for $E(X-a)_+$ when $\alpha > 1$.
\end{abstract}


\end{frontmatter}

\section{Introduction}

A univariate stable r.v. $Z$ with index $\alpha$,
skewness $\beta$, scale $\gamma$, and location $\delta$ has characteristic function
\begin{equation}\label{eq:stable.cf.1.param}
\phi(u) = \phi(u|\alpha,\beta,\gamma,\delta)= E \exp(i u Z) = \exp(-\gamma^\alpha [ |u|^\alpha + i \beta \eta(u,\alpha)] + i u \delta ),
\end{equation}
where $0 < \alpha \le 2$, $-1 \le \beta \le 1$, $\gamma >0$, $\delta \in \reals$ and
$$\eta(u,\alpha) =  \cases{
  - (\sign u) \tan (\pi \alpha/2) |u|^\alpha  & $\alpha \ne 1$ \cr
  (2/\pi) \, u \ln |u| & $\alpha=1$.}$$
In the notation of \citetwo{samorodnitsky:taqqu:1994}{Samorodnitsky, G.}{Taqqu, M.}, this is a
$S_\alpha(\gamma,\beta,\delta)$ distribution.  We will use the notation
$X \sim \Stable \alpha \beta {,\gamma,\delta;1}$ (the ``;1'' is used to distinguish
between this parameterization and a continuous one used below).


The purpose of this paper is to derive expressions for truncated fractional moments $E X_+^p = E (X \indicator {{ \{X \ge 0 \} }} )^p$ for general
stable laws.  To do this, define the functions for real $x$ and  $d$
\begin{eqnarray*}
g_d(x|\alpha,\beta) & = &   \cases{
   \displaystyle{ \int_0^\infty \cos(x r + \beta \eta(r,\alpha)) \, r^{d-1} e^{-r^\alpha} dr } & $ 0 < d < \infty$ \cr
  \displaystyle{\int_0^\infty [\cos(x r + \beta \eta(r,\alpha)) - 1] \,  r^{d-1}  \, e^{-r^\alpha} dr } & $ -2\min(1,\alpha) < d \le 0$ } \\
\gt_d(x|\alpha,\beta) & = & \cases{
    \displaystyle{\int_0^\infty \sin(x r + \beta \eta(r,\alpha)) \, r^{d-1} e^{-r^\alpha} dr } & $ -\min(1,\alpha) < d < \infty $\cr
    \displaystyle{\int_0^\infty [\sin(x r + \beta \eta(r,\alpha)) -xr ] \, r^{d-1} e^{-r^\alpha} dr } & $\alpha >1$, $ -\alpha < d \le -1$. }
\end{eqnarray*}
The functions $g_d(\cdot|\alpha,\beta)$ and $\gt_d(\cdot|\alpha,\beta)$, for integer subscripts $d=1,2,3,\ldots$ were introduced in
\citetwo{abdulhamid:nolan:1998}{Abdul-Hamid, H.}{Nolan, J. P.}.
(The notation was slightly different there: a factor of $(2 \pi)^{-d}$ was included in the definition
and $g_{\alpha,d}(x,\beta)$ was used instead of $g_d(x|\alpha,\beta)$,
while $q_{\alpha,1}(x,\beta)$ was used instead of $\gt_1(x|\alpha,\beta)$.)

The expressions for $EX_+^p$  will involve the functions $g_{-p}(\cdot|\alpha,\beta)$ and $\gt_{-p}(\cdot|\alpha,\beta)$, i.e. negative values fractional values of the subscript $d$.
Before proving that result, we show that the functions $g_d(\cdot|\alpha,\beta)$ and $\gt_d(\cdot|\alpha,\beta)$ have multiple uses.
For a standardized univariate stable law,  Fourier inversion of the characteristic function shows that the
d.f. and density are given by
\begin{eqnarray}
F(x|\alpha,\beta)-F(0|\alpha,\beta) & = & \frac{1}{\pi} \left( \gt_0(x|\alpha,\beta) - \gt_0(0|\alpha,\beta) \right) \label{eq:univar.cdf.by.gt0} \\
f(x|\alpha,\beta) & = &  {1 \over \pi}   g_1(x|\alpha,\beta). \nonumber
\end{eqnarray}
We note that there are explicit formulas for $F(0|\alpha,\beta)$
when $\alpha \ne 1$.

The $g_d(\cdot|\alpha,\beta)$ functions are used in a similar way to give $d$-dimensional
stable densities, see  Theorem~1 of \citetwo{abdulhamid:nolan:1998}{Abdul-hamid, H.}{Nolan, J. P.}
(note that there is a sign mistake in that formula when $\alpha=1$),
and \mycite{nolan:2017a}{Nolan, J. P.}  uses both $g_d(\cdot|\alpha,\beta)$ and $\gt_d(\cdot|\alpha,\beta)$  to give an expression for
multivariate stable probabilities.

Another use of these functions is in conditional expectation of $X_2$ given $X_1=x$ when $(X_1,X_2)$ are jointly stable with
zero shift and spectral measure $\Lambda$.  In general,
the conditional expectation is a complicated non-linear function; here it is restated in terms of these functions.
If $\alpha > 1$ or ($\alpha \le 1$ and (5.2.4) in \citetwo{samorodnitsky:taqqu:1994}{Samorodnitsky, G.}{Taqqu, M.} holds), then
Theorems 5.2.2 and 5.2.3 in \citetwo{samorodnitsky:taqqu:1994}{Samorodnitsky, G.}{Taqqu, M.} show that the conditional expectation
exists for $x$ in the support of $X_1$ and is given by
\begin{eqnarray*}
\lefteqn{ E(X_2|X_1=x) = } \\
& & \cases{
c_1 x + c_2 \left[ \displaystyle{\frac{1 - (x/\gamma_1) \gt_1(x/\gamma_1|\alpha,\beta_1)}{g_1(x/\gamma_1|\alpha,\beta_1)/\gamma_1}}\right] & $\alpha \ne 1$ \cr
c_0 + c_1 \displaystyle{ \frac{x-\mu_1}{\gamma_1}} + c_2 \, \displaystyle{ \frac{\gt_1((x-\mu_1)/\gamma_1 - (2 \beta_1/\pi ) \ln \gamma_1 |1,\beta_1)}
                  {g_1(x/\gamma_1|1,\beta_1)}} & $\alpha=1, \beta_1 \ne 0$ \cr
c_0 + c_1 \displaystyle{ \frac{x-\mu_1}{\gamma_1}} & $ $ \cr
           \quad + c_2 \,  \left[ \displaystyle{ \frac{(1-\ln \gamma_1) g_1((x-\mu_1)/\gamma_1|1,0)
           +h_1((x-\mu_1)/\gamma_1|1,0)}{ g_1(x/\gamma_1|1,0)}} \right]  & $\alpha=1, \beta_1 = 0$, }
\end{eqnarray*}
where $\beta_1$ and $\gamma_1$  are the skewness and scale parameters of $X_1$,
 and the constants and function $h_1(\cdot|1,0)$ are given by
\begin{eqnarray*}
c_0 & = & -\frac{2}{\pi} \int_\sphere s_2 \ln |s_1| \, \Lambda(d\bolds) \\
c_1 & = & \cases{
      \displaystyle{ \frac{ \kappa_1 +\beta_1 \tan^2 (\pi \alpha/2) \kappa_2 } {\gamma_1^{\alpha} (1+\beta_1^2 \tan^2  (\pi \alpha/2) )  }} & $\alpha \ne 1$ \cr
      \kappa_2/ \beta_1 & $\alpha=1, \beta_1 \ne 0$ \cr
      \kappa_1 & $\alpha=1,\beta_1=0$} \\
c_2 & = & \cases{
     \displaystyle{\frac{\tan (\pi \alpha/2) (\kappa_2 - \beta_1 \kappa_1)}{\gamma_1^{\alpha} (1+\beta_1^2 \tan^2  (\pi \alpha/2) )  }} & $\alpha \ne 1$ \cr
      (\kappa_2-\beta_1 \kappa_1)/ \beta_1 & $\alpha=1, \beta_1 \ne 0$ \cr
      - 2 \kappa_2/\pi  & $\alpha=1,\beta_1=0$} \\
\kappa_1 & = & \covariation {X_2}{X_1}{\alpha}  =  \cases{
        \int_{\sphere} s_2 \signedpower{s_1}{\alpha-1} \, \Lambda(d\bolds) & $ \alpha \ne 1 $ \cr
        \int_{\sphere} s_2 \signedpower{s_1}{0} \, \Lambda(d\bolds) = \int_{\sphere} s_2 \, \mathrm{sign}(s_1) \, \Lambda(d\bolds) & $ \alpha = 1$ } \\
\kappa_2 & = & \int_{\sphere} s_2 |s_1|^{\alpha-1} \, \Lambda(d\bolds) \\
\mu_1 & = &  -\frac{2}{\pi} \int_\sphere s_1 \ln |s_1| \, \Lambda(d\bolds) \\
h(x|1,0) & = &  \int_0^\infty \cos(xr) (\log r)  e^{-r} dr.
\end{eqnarray*}
In the terms above, $\sphere$ is the unit circle and $ \covariation {X_2}{X_1}{\alpha}$ is the $\alpha-$covariation.
Note that if $\Lambda$ is symmetric, then $c_0=\kappa_2=\beta_1=\mu_1=0$, so  $c_2=0$ and
$$E(X_2|X_1=x) = \frac{ \covariation {X_2}{X_1}{\alpha}}{\gamma_1^\alpha} \,  x$$
is linear.

\section{Truncated moments $E X_+^p$}

The main result of this paper is the following expression for the fractional truncated moment of a stable r.v.

\startthm\label{thm:moment.of.truncation}
Let $X \sim \Stable \alpha \beta {,\gamma,\delta;1}$ with any $0 < \alpha < 2$ and any $-1 \le \beta \le 1$ and set
$$\delta^* = \cases{ \delta/\gamma & $\alpha \ne 1$ \cr
\delta/\gamma + (2/\pi)\beta \log \gamma & $\alpha = 1$.} $$
For $ p < \alpha$, define $m^p(\alpha,\beta,\gamma,\delta) = EX_+^p$.  \\

\myitem{(a)}
When $p=0$,
$$m^0(\alpha,\beta,\gamma,\delta)= P(X > 0) = \frac{1}{2} - \frac 1 \pi \gt_0(-\delta^*|\alpha,\beta).$$
When $0 < p < \min(1,\alpha)$,
\begin{eqnarray*}
m^p(\alpha,\beta,\gamma,\delta)& = & \gamma^p \frac{\Gamma(p+1)}{\pi} \left[   \sin\parenspipbytwo \left( \frac{\Gamma(1-p/\alpha)}{p}
        - g_{-p}(-\delta^*|\alpha,\beta) \right)  \strutfive \right. \\
 & &  \hspace{4cm} \left. - \cos\parenspipbytwo \gt_{-p}(-\delta^*|\alpha,\beta) \strutfive \right] .
\end{eqnarray*}
When $p=1 < \alpha < 2$,
$$ m^p(\alpha,\beta,\gamma,\delta) = \gamma  \left[  \frac{\delta^*} 2 + \frac 1 \pi   \left(\Gamma(1-1/\alpha) - g_{-1}(-\delta^*|\alpha,\beta) \right) \right] .$$
When $1 < p < \alpha < 2$,
\begin{eqnarray*}
m^p(\alpha,\beta,\gamma,\delta) & = & \gamma^p \frac{\Gamma(p+1)}{\pi} \left[ \sin\parenspipbytwo \left( \frac{\Gamma(1-p/\alpha)}{p} - g_{-p}(-\delta^*|\alpha,\beta)\right) \strutfive \right. \\
  & & \hspace{3cm} \left.  +  \cos\parenspipbytwo
          \left(  \frac{\delta^*} {\alpha} \Gamma((1-p)/\alpha) - \gt_{-p}(-\delta^*|\alpha,\beta) \right) \strutfive \right] .
\end{eqnarray*}
\myitem{(b)}  $E X_-^p = E (-X)_+^p = m^p(\alpha,-\beta,\gamma,-\delta).$
\stopthm

\startpf {} \myitem{(a)} To simplify calculations, first assume $\gamma=1$; the adjustment for $\gamma \ne 1$ is discussed below.
 When $p=0$, $EX_+^0 = \int_0^\infty 1 \, f(x) dx = P(X>0)$, and (\ref{eq:univar.cdf.by.gt0})
and $\gt_0(x|\alpha,\beta) \to \pi/2$ as $x \to \infty$ gives the value in terms of $\gt_0(\cdot|\alpha,\beta)$. \\
When $0 < p < \min(1,\alpha)$, Corollary~2 of \mycite{pinelis:2011}{Pinelis, I.} with $k=\ell=0$ shows
\begin{equation}\label{eq:frac.momment.p.lt.1}
E X_+^p = \frac{\Gamma(p+1)}{\pi} \int_0^\infty \Re \frac{\phi(u)-1}{(iu)^{p+1}} du.
\end{equation}
First assume  $\alpha \ne 1$ and set  $\zeta=\zeta(\alpha,\beta) = -\beta \tanconst$ and restricting to $u > 0$,
\begin{eqnarray*}
\frac{\phi(u)-1}{(iu)^{p+1}} & = &  \left( \left[ e^{-u^\alpha (1+i\zeta )+i \delta u }-1 \right] (-i) e^{-i (\pi/2) p} \right) u^{-p-1}  \\
 & = &  \left( -i \left( e^{-u^\alpha} \left[ \cos(\delta u - \zeta u^\alpha) + i \sin(\delta u - \zeta u^\alpha)\right] -1  \right)
              e^{-i (\pi/2) p} \right) u^{-p-1}\\
 & = &  \left( \left[ e^{-u^\alpha} \sin(\delta u - \zeta u^\alpha) - i \left( e^{-u^\alpha} \cos(\delta u - \zeta u^\alpha)  -1 \right)  \right]
              \left[ \cos\parenspipbytwo - i \sin\parenspipbytwo \right] \right) u^{-p-1}
\end{eqnarray*}
And therefore
\begin{eqnarray*}
\Re \frac{\phi(u)-1}{(iu)^{p+1}}  & = & \left[  \cos\parenspipbytwo e^{-u^\alpha} \sin(\delta u - \zeta u^\alpha )
     - \sin\parenspipbytwo \left(  e^{-u^\alpha} \cos(\delta u - \zeta u^\alpha) -1 \right) \right] u^{-p-1} \\
     & = &  \cos\parenspipbytwo \sin(\delta u - \zeta u^\alpha ) u^{-p-1} e^{-u^\alpha} \\
     & &  \quad - \sin\parenspipbytwo \left( \left[ \cos(\delta u - \zeta u^\alpha) -1 \right] u^{-p-1}  e^{-u^\alpha} +  (e^{-u^\alpha}-1) u^{-p-1} \right)
\end{eqnarray*}
Integrating this from 0 to $\infty$, substituting $t=u^\alpha$ in the last term to get
$$E X_+^p = \frac{\Gamma(p+1)}{\pi} \left[ -\cos\parenspipbytwo \gt_{-p}(-\delta|\alpha,\beta) - \sin\parenspipbytwo
    \left\{ g_{-p}(-\delta|\alpha,\beta) - \Gamma(1-p/\alpha)/p \right\} \right].$$

Next consider $0 < p < \alpha = 1$.  Use (\ref{eq:frac.momment.p.lt.1}) again, so we need to simplify
\begin{eqnarray*}
\frac{\phi(u)-1}{(iu)^{p+1}} & = &  \left( \left[ e^{-u (1+i \beta \eta(u,1))+i \delta u }-1 \right] (-i) e^{-i (\pi/2) p} \right) u^{-p-1}  \\
 & = &  \left( -i \left( e^{-u} \left[ \cos(\delta u - \beta \eta(u,1) ) + i \sin(\delta u - \beta \eta(u,1) )\right] -1  \right)
              e^{-i (\pi/2) p} \right) u^{-p-1}\\
 & = &  \left( \left[ e^{-u} \sin(\delta u - \beta \eta(u,1)) - i \left( e^{-u} \cos(\delta u - \beta \eta(u,1))  -1 \right)  \right] \right. \\
 & & \hspace{1cm} \times   \left. \left[ \cos\parenspipbytwo - i \sin\parenspipbytwo \right] \right) u^{-p-1} \\
 \Re \frac{\phi(u)-1}{(iu)^{p+1}} & = &   \left[ \cos\parenspipbytwo  e^{-u} \sin(\delta u - \beta \eta(u,1)) \strutfive \right. \\
     &  & \hspace{1cm} \left. - \sin\parenspipbytwo \left( e^{-u} ( \cos(\delta u - \beta \eta(u,1))  -1) + (e^u -1 ) \right) \strutfive \right]  u^{-p-1}.
\end{eqnarray*}
Integrating from $0$ to $\infty$ yields
$$E X_+^p = \frac{\Gamma(p+1)}{\pi} \left[ -\cos\parenspipbytwo \gt_{-p}(-\delta|1,\beta) - \sin\parenspipbytwo
    \left\{ g_{-p}(-\delta|1,\beta) - \Gamma(1-p)/p \right\} \right].$$

When $p= 1 < \alpha < 2$, $E X$ exists and is equal to $\delta$.  Using Corollary~2 of \mycite{pinelis:2011}{Pinelis, I.} with $k=1$, $\ell=0$  shows
$$ E X_+  = \frac{1}{2} EX + \frac{\Gamma(2)}{\pi} \int_0^\infty \Re \frac{\phi(u)-1}{(iu)^{p+1}} du
     = \frac{\delta}{2} + \frac{1}{\pi} \int_0^\infty \Re \frac{\phi(u)-1}{(iu)^{p+1}} du.$$
The integrand is the same as above, with $\cos \parenspipbytwo=0$ and $\sin\parenspipbytwo=1$, so
$$ E X_+ =  \frac{\delta}{2} - \frac{1}{\pi} \left[ g_{-1}(-\delta|\alpha,\beta) - \Gamma(1-1/\alpha) \right].$$

When $1 < p < \alpha < 2$, Corollary~2 of \mycite{pinelis:2011}{Pinelis, I.} with $k=\ell=1$  shows
\begin{equation}\label{eq:trunc.moment.p.gt.1}
E X_+^p = \frac{\Gamma(p+1)}{\pi} \int_0^\infty \Re \frac{\phi(u)-1-i u EX}{(iu)^{p+1}} du,
\end{equation}
Since $\alpha > 1$, $EX$ exists and is equal to $\delta$.  As above, for $u > 0$,
\begin{eqnarray*}
\frac{\phi(u)-1-i u \delta}{(iu)^{p+1}} & = &  \left( \left[ e^{-u^\alpha (1+i\zeta u^\alpha)+i \delta u }-1-i \delta u \right] (-i) e^{-i (\pi/2) p} \right) u^{-p-1}  \\
 & = &  \left( -i \left( e^{-u^\alpha} \left[ \cos(\delta u - \zeta u^\alpha) + i \sin(\delta u - \zeta u^\alpha)\right] -1 -i \delta u \right)
              e^{-i (\pi/2) p} \right) u^{-p-1}\\
 & = &   \left[ \left( e^{-u^\alpha} \sin(\delta u - \zeta u^\alpha) - \delta u \right) - i \left( e^{-u^\alpha}  \cos(\delta u - \zeta u^\alpha)  -1 \right)  \right]\\
 & &      \quad \times \left[ \cos\parenspipbytwo - i \sin\parenspipbytwo \right]  u^{-p-1}
\end{eqnarray*}
And therefore
\begin{eqnarray*}
\Re \frac{\phi(u)-1-iu \delta}{(iu)^{p+1}}
     & = &  \left[  \cos\parenspipbytwo\left(  e^{-u^\alpha} \sin(\delta u - \zeta u^\alpha )-\delta u \right) \right. \\
     & & \quad    \left.  - \sin\parenspipbytwo \left(  e^{-u^\alpha} \cos(\delta u - \zeta u^\alpha) -1 \right) \right] u^{-p-1} \\
     & = &  \cos\parenspipbytwo \left( [ \sin(\delta u - \zeta u^\alpha ) - \delta u ] u^{-p-1} e^{-u^\alpha}  + \delta (e^{-u^\alpha}-1) u^{-p}   \right) \\
     & &  \quad - \sin\parenspipbytwo \left( \left[ \cos(\delta u - \zeta u^\alpha) -1 \right] u^{-p-1}  e^{-u^\alpha} +  (e^{-u^\alpha}-1) u^{-p-1} \right)
\end{eqnarray*}
Plugging this into (\ref{eq:trunc.moment.p.gt.1}) and integrating yields
\begin{eqnarray*}
E X_+^p
     & = & \frac{\Gamma(p+1)}{\pi} \left\{  \cos\parenspipbytwo \left[ -\gt_{-p}(-\delta|\alpha,\beta))
                 + (\delta/\alpha) (\Gamma((1-p)/\alpha)  \right] \strutfive \right. \\
     & &  \hspace{2cm} \left. + \sin\parenspipbytwo \left[ -g_{-p}(-\delta|\alpha,\beta)  + \Gamma(1-p/\alpha)/p \right] \strutfive \right\}
\end{eqnarray*}

Now consider $\gamma \ne 1$.  If $X \sim \Stable \alpha \beta {,\gamma, \delta;1}$, then $X \eqd \gamma Y$, where
$Y \sim \Stable \alpha \beta {,1,\delta^*;1}$, so $E X_+^p = \gamma^p E Y_+^p$. In symbols,
$$m^p(\alpha,\beta,\gamma,\delta) = \gamma^p m^p(\alpha,\beta,1,\delta^*).$$

\myitem{(b)} This follows from $-X \sim \Stable \alpha {-\beta} {,\gamma,-\delta;1} $.
\stoppf

When $-1 < p < 0$, we conjecture that
$$m^p(\alpha,\beta,\gamma,\delta) = \gamma^p \frac{\Gamma(p+1)}{\pi} \left[   - \sin \parenspipbytwo g_{-p}(-\delta^*|\alpha,\beta)
   - \cos\parenspipbytwo \gt_{-p}(-\delta^*|\alpha,\beta) \right] . $$

\section{Related results}

There are several corollaries to the preceding result.  First, taking $p=1$ in the previous result shows the following.
\startcor If $X \sim \Stable \alpha \beta {,\gamma,\delta;1}$ with $\alpha > 1$, $-1 \le \beta \le 1$, $a \in \reals$
$$E ( X-a)_+ = \frac {\delta-a} 2 + \frac{\gamma}{\pi} \left[  \Gamma \left(1-\frac 1 \alpha \right) - g_{-1} \left.  \left(\frac{\delta-a} \gamma \right| \alpha,\beta \right)   \strutfive \right]. $$
\stopcor

\comment{ **********************************
\bigskip
OLD \startcor\label{cor:abs.sgn.momentsOLD} If $X \sim \Stable \alpha \beta {,\gamma,\delta;1}$ with $0 <\alpha < 2$, $-1 \le \beta \le 1$, $0 \le p < \alpha$.
\begin{eqnarray*}
E|X|^p & = & \cases{
    \displaystyle{ -\gamma^p K   g_{-p}(-\delta^*|\alpha,\beta) }  & $-1 < p < 0$ \cr
    1 & $p=0$ \cr
    \displaystyle{ \gamma^p K   \left( \frac{\Gamma(1-p/\alpha)}{p} - g_{-p}(-\delta^*|\alpha,\beta) \right)} & $ 0 < p < \alpha$
    } \\
E X\sgnpower{p} & = & \cases{
     \displaystyle{  -\frac{2}{\pi} \gt_0(-\delta^*|\alpha,\beta) } & $p=0$ \cr
     \displaystyle{ - \gamma^p \widetilde{K}
              \gt_{-p}(-\delta^*|\alpha,\beta)  }  & $p \in (-1,0) \cup (0,\min(1,\alpha))$ \cr
     \delta & $p=1$ \cr
     \displaystyle{ \gamma^p \widetilde{K}
              \left(  \frac{\delta^* \Gamma( (1-p)/\alpha)}{\alpha} - \gt_{-p}(-\delta^*|\alpha,\beta) \right)} & $ 1 < p < \alpha$,
  }
\end{eqnarray*}
where $K=K(p) = \frac{2 \Gamma(p+1)}{\pi} \sin\parenspipbytwo$ and $\widetilde{K}=\widetilde{K}(p) = \frac{2 \Gamma(p+1)}{\pi} \cos\parenspipbytwo$.
\stopcor
  ********** END COMMENT }

Combining parts (a) and (b) of Theorem~\ref{thm:moment.of.truncation} yields.
\startcor\label{cor:abs.sgn.moments} If $X \sim \Stable \alpha \beta {,\gamma,\delta;1}$ with $0 <\alpha < 2$, $-1 \le \beta \le 1$, $-1 \le p < \alpha$.
\begin{eqnarray*}
E|X|^p & = & \displaystyle{ \gamma^p \frac{2 \Gamma(p+1)}{\pi} \sin\parenspipbytwo  \left(  \frac{\delta^* \Gamma( 1-p/\alpha)}{p} \indicator{ {\{ p > 0 \}} }  - g_{-p}(-\delta^*|\alpha,\beta) \right)}  \\
E X\sgnpower{p} & = & \displaystyle{ \gamma^p \frac{2 \Gamma(p+1)}{\pi} \cos\parenspipbytwo
              \left(  \frac{\delta^* \Gamma( (1-p)/\alpha)}{\alpha} \indicator{ { \{p > 1} \}} - \gt_{-p}(-\delta^*|\alpha,\beta) \right). }
\end{eqnarray*}
\stopcor

\startpf {}
$E|X|^p = E X_-^p + E X_+^p = m^p(\alpha,\beta,\gamma,\delta) + m^p(\alpha,-\beta,\gamma,-\delta)$ and
$E X\sgnpower{p} = m^p(\alpha,\beta,\gamma,\delta) -  m^p(\alpha,-\beta,\gamma,-\delta)$.  Use Theorem~\ref{thm:moment.of.truncation}
and the reflection property: $g_d(-x|\alpha,\beta)=g_d(x|\alpha,-\beta)$.  Note that as $p \to 0$, $E|X|^p \to E \, 1 = 1$ and
$E X\sgnpower{p} \to -(2/\pi) \gt_{0}(\delta^*|\alpha,\beta)=P(X>0)-P(X<0)=1-2 F(0)$.
Also as $p \to 1$, $E X\sgnpower{p} \to \delta$.
\stoppf

\comment{  *****************************************************************************************************

\startcor  If $X \sim \Stable \alpha \beta {,\gamma,\delta;1}$ with $0 <\alpha < 2$, $-1 \le \beta \le 1$.
Then the moment generating functions of $\log |X|$, $\log X_+$ and $\log X_-$ for $u \ge 0$ are
\begin{eqnarray*}
M_+(u) & = & E \exp( u \log X_+ ) = m^u(\alpha,\beta,\gamma,\delta) \\
M_-(u) & = &  E \exp( u \log X_- ) = m^u(\alpha,-\beta,\gamma,-\delta) \\
M(u)   & = & E \exp( u \log |X| ) = m^u(\alpha,\beta,\gamma,\delta) + m^u(\alpha,-\beta,\gamma,-\delta).
\end{eqnarray*}
The means and variances are
$$\begin{array}{llllll}
E \log X_+ & = & ???  \hspace{1cm} \mathrm{Var}(\log X_+) & = & ??? \\
E \log X_- & = & ???  \hspace{1cm} \mathrm{Var}(\log X_-) & = & ??? \\
E \log |X| & = & ???  \hspace{1cm} \mathrm{Var}(\log |X|) & = & ??? \\
\end{array}$$
\stopcor
\startpf {} $E \log X_+ = M_+^\prime(0) = ...$.  \unfinished

   ***************************************************************************************************** }

In the strictly stable case, the expressions for $EX_+^p$ can be simplified using
closed form expressions for $g_d(0|\alpha,\beta)$ and $\gt_d(0|\alpha,\beta)$ when $\alpha \ne 1$.
To state the result, set
$$\theta_0 = \theta_0(\alpha,\beta) = \cases{
    \alpha^{-1} \arctan\left(\beta \tanconst\right) & $\alpha \ne 1$ \cr
    \pi/2 & $\alpha=1$.} $$

\startlem \label{lem:g0.gt0} When $\alpha \ne 1$,
\begin{eqnarray*}
  g_d(0|\alpha,\beta) &=& \cases{
       (\cos \alpha \theta_0)^{d/\alpha}  \cos (d \, \theta_0) \Gamma(1+d/\alpha)/d & $d > 0$\cr
       (\ln (\cos \alpha \theta_0))/\alpha & $d=0$ \cr
       \left[ (\cos \alpha \theta_0)^{d/\alpha}  \cos (d \, \theta_0)-1 \right] \Gamma(1+d/\alpha)/d  & $-\alpha < d < 0$  } \\
\gt_d(0|\alpha,\beta) &=& \cases{
   - (\cos \alpha \theta_0)^{d/\alpha}  \sin (d \, \theta_0) \Gamma(1+d/\alpha)/d & $ d \in  (-\alpha,0) \cup (0,\infty)$  \cr
   - \theta_0 & $d=0$. }
\end{eqnarray*}
\stoplem

\startpf {} Substitute $u=r^\alpha$ in the expressions for $g_d(0|\alpha,\beta)$ and $\gt_d(0|\alpha,\beta)$.
Then use respectively the integrals 3.944.6, 3.948.2, 3.945.1, 3.944.5, and 3.948.1  pg. 492-493 of
\citetwo{gradshteyn:ryzhik:2000}{Gradshteyn, I.}{Ryzhik, I.}. (Note that some of these formulas have
mistyped exponents.) Finally, when $\alpha \ne 1$, $\alpha \theta_0 = - \arctan \zeta$,
and for the allowable values of $\alpha$ and $\theta_0$,
$$\cos \alpha \theta_0 = |\cos \alpha \theta_0|
 =  (1+\tan^2 \alpha \theta_0)^{-1/2} = (1+\zeta^2)^{-1/2}.$$
\stoppf

The following is a different proof of Theorem~2.6.3 of \mycite{zolotarev:1986}{Zolotarev, V.}.  

\startcor \label{cor:moment.of.truncation.strictly.stable} Let $X$ be strictly stable, e.g. $X \sim \Stable \alpha \beta {,\gamma,0;1}$ with $\alpha \ne 1$
 or ($\alpha=1$ and  $\beta=0$) and  $0 < p < \alpha$. \\
\myitem{(a)} The fractional moment of the positive part of $X$ is
$$EX_+^p =   \frac{\gamma^p }{ (\cos \alpha \theta_0)^{p/\alpha}}
{\Gamma(1-p/\alpha) \over  \Gamma(1-p) } { \sin p(\pi/2 +\theta_0)
        \over  \sin (p \pi)}.$$
\myitem{(b)} The fractional moment of the negative part of $X$ is $E X_-^p = E (-X)_+^p$,
which can be obtained from the right hand side above by replacing $\theta_0$ with $-\theta_0$.
\stopcor
When $p=1$, the product  $\Gamma(1-p) \sin (\pi p)$ in the denominator above  is interpreted as the limiting value as $p \to 1$, which is $\pi$. \\

\startpf {} Note that when $X$ is strictly stable, $\delta^*=0$.  First assume $0 < p < \min(1,\alpha)$ and substitute  Lemma~\ref{lem:g0.gt0} into this case of  Theorem~\ref{thm:moment.of.truncation}
\begin{eqnarray*}
E X_+^p & = & \frac{\gamma^p \Gamma(p+1)}{\pi}  \left[ \sin(\pi p/2) \left( \frac{\Gamma(1-p/\alpha)}{p}
             - \left( (\cos \alpha \theta_0)^{-p/\alpha} \cos(-p \theta_0) - 1\right)\frac{\Gamma(1-p/\alpha)}{-p}    \right) \right.  \\
      & & \hspace{2cm} \left. - \cos(\pi p/2) (-\cos (\alpha \theta_0)^{-p/\alpha} \sin (-p \theta_0) \frac{\Gamma(1-p/\alpha)}{-p}      \right] \\
     & = & \frac{\gamma^p \Gamma(p+1)\Gamma(1-p/\alpha) }{\pi p (\cos \alpha \theta_0)^{p/\alpha} }  \left[ \sin(\pi p/2) \cos(p \theta_0)
                   + \cos(\pi p/2) \sin (p \theta_0 \right] \\
     & = & \frac{\gamma^p \Gamma(p+1)\Gamma(1-p/\alpha) }{\pi p (\cos \alpha \theta_0)^{p/\alpha} }  \sin(\pi p/2+p \theta_0)
\end{eqnarray*}
Using the identity $\Gamma(p+1)= \pi p/(\Gamma(1-p) \sin p \pi) $ gives the result.

When $p=1 < \alpha$, again using the appropriate part of Theorem~\ref{thm:moment.of.truncation} shows
\begin{eqnarray*}
E X_+ & = & \gamma   \left[ 0 + \frac{1}{\pi} \left( \Gamma(1-1/\alpha) -  \left( (\cos \alpha \theta_0)^{-1/\alpha} \cos(- \theta_0) - 1 \right)\frac{\Gamma(1-1/\alpha)}{-1}  \right) \right]  \\
     & = & \frac{ \gamma \Gamma(1-1/\alpha)}{\pi} \left[ \cos \alpha \theta_0)^{-1/\alpha} \cos(\theta_0) \right].
\end{eqnarray*}
When $1 < p < \alpha$, using Theorem~\ref{thm:moment.of.truncation} and $\delta^*=0$,
\begin{eqnarray*}
E X_+^p & = & \frac{\gamma^p \Gamma(p+1)}{\pi}  \left[ \sin(\pi p/2) \left( \frac{\Gamma(1-p/\alpha)}{p}
             - \left((\cos \alpha \theta_0)^{-p/\alpha} \cos(-p \theta_0) - 1\right)\frac{\Gamma(1-p/\alpha)}{-p}    \right) \right.  \\
      & & \hspace{2cm} \left. + \cos(\pi p/2)  \left(0-\cos (\alpha \theta_0)^{-p/\alpha} \sin (-p \theta_0) \right) \frac{\Gamma(1-p/\alpha)}{-p}      \right],
\end{eqnarray*}
and the rest is like the first case.
\stoppf

The standard parameterization used above is discontinuous in the parameters near $\alpha=1$, and it is not a scale-location family 
when $\alpha=1$.  To avoid this, a continuous parameterization that is a scale-location family can be used.
We will say $X \sim \Stable \alpha \beta {,\gamma,\delta;0}$ if it has characteristic function 
$$ E \exp(i u X ) =  \cases{
\exp \left(
  -\gamma^\alpha |u|^\alpha \left[1 + i \beta (\tanconst) (\sign u)
  (|\gamma u|^{1-\alpha}-1) \right] + i \delta u  \right)
   & $\alpha \ne 1$ \cr
\exp \left(
  -\gamma |u|        \left[ 1 + i \beta (2/\pi) (\sign u) \log(\gamma |u|)
    \right] + i \delta u  \right)           & $\alpha=1$.}$$
A stable r. v. $X$ can be expressed in both the 0-parameterization and the 1-parameterization, in which case the 
index $\alpha$, the skewness $\beta$ and the scale $\gamma$ are the same. The only difference is
in the location parameter: if $X$ is simultaneously $\Stable \alpha \beta {,\gamma,\delta_0;0}$ and 
$\Stable \alpha \beta {,\gamma,\delta_1;1}$, then the shift parameters are related by 
$$\delta_1 = \cases{ \delta_0 - \beta \gamma \tanconst & $\alpha \ne 1$ \cr     
                     \delta_0 - (2/\pi) \beta \gamma \log \gamma  & $\alpha = 1$.}$$
Therefore, if $X \sim \Stable \alpha \beta {,\gamma,\delta_0;0}$, 
$$E X_+^p = \cases{ m^p(\alpha,\beta,\gamma,\delta_0- \beta \gamma \tanconst  ) & $\alpha \ne 1$ \cr     
                    m^p(\alpha,\beta,\gamma,\delta_0 - (2/\pi) \beta \gamma \log \gamma  ) & $\alpha=1$.}$$
This quantity is continuous in all parameters.

For the above expressions for $EX_+^p$ to be of practical use, one must be able to evaluate $g_d(\cdot|\alpha,\beta)$ and 
$\gt_d(\cdot|\alpha,\beta)$.  When $d$ is a nonnegative integer, \mycite{nolan:2017a}{Nolan, J.} gives Zolotarev type integral expressions for these functions.
However, this is not helpful here, where negative, non-integer values of $d$ are needed.  We have written a short R program to 
numerically evaluate the defining integrals for $g_d(\cdot|\alpha,\beta)$ and $\gt_d(\cdot|\alpha,\beta)$.  A single evaluation takes 
less than 0.0002 seconds on a modern desktop. This faster than numerically evaluating $E X_+^p = \int_{0}^{\infty} x^p f(x|\alpha,\beta,\gamma,\delta) dx$,
because the latter requires many numerical calculations of the density $f(x|\alpha,\beta,\gamma,\delta)$.  

\section*{References}

\bibliography{TruncatedFractionalMoments}

\end{document}